\documentstyle[amssymb]{article}
\date{}
\begin{document}
{\noindent\Large
\texttt{A property of cyclotomic polynomials}
\footnote{AMS MCS 11C08.}}

\bigskip\noindent
{{Giovanni Falcone}\footnote{Supported by M.I.U.R., Universit\`a di Palermo.}\\
\footnotesize{Dipartimento di Metodi e Modelli Matematici}\vspace*{-1,5mm}\\
\footnotesize{Viale delle Scienze Ed. 8, I-90128 Palermo (Italy)}\vspace*{-1,5mm}\\
\footnotesize{gfalcone@unipa.it}}

\vskip1cm\begin{abstract}\noindent
Given two cyclotomic polynomials $\Phi_n(x)$ and $\Phi_m(x)$,
$n\not= m$, we determine the minimal natural number $k$ such that
we can write
$$k=a(x)\Phi_n (x)+b(x)\Phi_m(x),$$
with $a(x)$ and $b(x)$ integer polynomials.\end{abstract}

\large\vskip1cm\noindent In June 2000 the author of this note gave
a talk at the conference "Combinatorics 2000" in Gaeta (Italy),
titled \em Dividing Cyclotomic Polynomials\em . The content of
that talk was proposed to the {\bf American Mathematical Monthly}
as Problem 10914 in Volume 109, January 2002, and a composite
solution by R.~Stong and N.~Komanda was published on Volume 110,
October 2003. The problem was also solved by R.~Chapman,
C.~P.~Rupert, the GCHQ Problem Solving Group and the proposer,
whose original solution is given in this note.

\vskip1cm\noindent{\bf A property of cyclotomic polynomials.}
\bigskip
We denote, as usual, with $\Phi_k(x)$ the cyclotomic polynomial
(with integer coefficients) of index $k$, defined inductively
through the identity
$$x^m-1=\prod \Phi_k(x),$$
where $k$ runs among the divisors of $m$.

\smallskip
The basic properties of cyclotomic polynomials are well known, as
well as the role they play in several branches of Mathematics.
Here we just recall that $\Phi_m(1)=p$, if $m$ is a $p$-power ($p$
prime), $\Phi_m(1)=1$ otherwise.

\smallskip
Given two different cyclotomic polynomials $\Phi_m(x),\Phi_n(x)$,
we can find two polynomials $s(x),t(x)\in{\Bbb Q}[x]$ such that
$$1=s(x)\Phi_m(x)+t(x)\Phi_n(x),$$
$\Phi_m(x)$ and $\Phi_n(x)$ being irreducible as element of the
euclidean ring ${\Bbb Q}[x]$. Since the ring ${\Bbb Z}[x]$ is a
unique factorization domain, but not a euclidean ring, here this
property fails. But the existence of an integer $k$ and of two
integer polynomials $a(x),b(x)\in {\Bbb Z}[x]$, such that
$$k=a(x)\Phi_m(x)+b(x)\Phi_n(x),$$
is manifestly guaranteed.

\bigskip
This note is essentially an integration of \cite{Finite}, where we
used cyclotomic polynomials to define minimal polynomials of
finite automorphisms of groups. The knowledge of the integer $k$
is useful in order to give a {\em unique} definition of the
minimal polynomial for elements of finite order in a ring, as well
as a decomposition of a group, on which a finite automorphism is
acting, into the direct sum of invariant subgroups. We refer to
\cite{Finite} for more details.

According to \cite{Finite}, the following proposition holds good:

\bigskip
\noindent{\bf Proposition.} {\em Let $m>n$ be two integers. If $n$
does not divide $m$ then two polynomial $a(x),b(x)\in{\Bbb Z}[x]$
exist, such that
$$1=a(x)\Phi_m(x)+b(x)\Phi_n(x).$$}

In order to evaluate $k$ we can confine ourself then to the case
in which $n$ is a divisor of $m$:

\bigskip
\noindent{\bf Proposition.} {\em Let $\Phi_m(x),\Phi_n(x)$ be two
cyclotomic polynomials, and let $n$ be a divisor of $m$. Then two
polynomials $a(x), b(x)\in {\Bbb Z}[x]$ exist, such that
$$k=a(x)\Phi_m(x)+b(x)\Phi_n(x),$$
where $k=1$, if $\frac{m}{n}$ is not a prime-power, or $k=p$, if
$\frac{m}{n}=p^t$ ($p$ a prime). In the latter case, $p$ is the
smallest positive integer with that property.}

\bigskip
\noindent {\it Proof.} Note that, since $n$ divides $m$, then
$\Phi_m(x)$ is a divisor of the polynomial
$\Phi_{\frac{m}{n}}(x^n)$. In fact, if $\omega\in {\Bbb C}$ is a
primitive $m$-th root of unity, then $\omega$ generates a cyclic
group of order $m$ and $\omega^n$ generates a cyclic group of
order $\frac{m}{n}$, that is $\omega^n$ is a primitive
$\frac{m}{n}$-th root of unity, hence
$$\Phi_{\frac{m}{n}}(\omega^n)=0.$$
Each root of $\Phi_m(x)$ is then a root of
$\Phi_{\frac{m}{n}}(x^n)$. Since $\Phi_m(x)$ is irreducible, this
means that $\Phi_m(x)$ is a divisor of $\Phi_{\frac{m}{n}}(x^n)$,
even over ${\Bbb Z}[x]$ by Gau\ss ' Lemma.

Write $\Phi_{\frac{m}{n}}(x^n)=a(x)\Phi_m(x)$ and let
$\sigma\in{\Bbb C}$ be a primitive $n$-th root of unity, so we
have
$$a(\sigma)\Phi_m(\sigma)=\Phi_{\frac{m}{n}}(\sigma^n)=\Phi_{\frac{m}{n}}(1)=k,
$$
where $k=1$, if $\frac{m}{n}$ is not a prime-power, or $k=p$, if
$\frac{m}{n}=p^t$ ($p$ prime).

Divide $a(x)\Phi_m(x)$ by $\Phi_n(x)$ and write
$$a(x)\Phi_m(x)=b(x)\Phi_n(x)+r(x),$$
where $r(x)$ is an integer polynomial of lower degree than the one
of $\Phi_n(x)$. We get
$$k=r(\sigma),$$
but this is possible only if $r(x)$ is a constant (equal to $k$),
since the degree of $r(x) - k$ is lower than the one of the
irreducible polynomial $\Phi_n(x)$.

In order to prove that, if $\frac{m}{n}=p^t$ ($p$ prime), no other
natural number $0<k<p$ is such that
$$k=a(x)\Phi_m(x)+b(x)\Phi_n(x),$$
we note that $\Phi_m(x)$ and $\Phi_n(x)$ have, in this case, the
same roots on a field {\sf K} of characteristic $p$ (see the
following Note), and for such a root $\alpha\in {\sf K}$ we would
have
$$k=a(\alpha)\Phi_m(\alpha)+b(\alpha)\Phi_n(\alpha)=0,$$
a contradiction. \hfill$\Box$

\bigskip
\noindent{\bf Note.} Let {\sf K} be a field of characteristic $p$,
$\alpha\in{\sf K}$ and $h=kp^r$, $p$ not dividing $k$. Then
$$\Phi_h(\alpha)=0\iff \Phi_k(\alpha)=0.$$
In fact, one can immediatily see that on ${\Bbb Z}[x]$ the
following two properties of cyclotomic polynomials hold good:

\noindent {\it i}) if $p$ is a prime, not dividing $n$, then
$\Phi_n(x)\Phi_{pn}(x)$ and $\Phi_n(x^p)$ have precisely the same
roots, hence
$$\Phi_{pn}(x)=\frac{\Phi_n(x^p)}{\Phi_n(x)};$$

\noindent {\it ii}) if $n=p_1^{r_1}\cdots p_s^{r_s}$ is the prime
factorization of $n$, then
$$\Phi_n(x)=
\Phi_{p_1\cdots p_s} (x^{ (p_1^{(r_1-1)}\cdots p_s^{(r_s-1)})
}).$$ In fact $\omega$ is a primitive $n$-th root of unity if and
only if $\omega^{(p_1^{(r_1-1)}\cdots p_s^{(r_s-1)})}$ is a
primitive $(p_1\cdots p_s)$-th root of unity.

\smallskip
Now let $k=p_1^{r_1}\cdots p_s^{r_s}$ be the prime factorization
of $k$, then $h=kp^r=p_1^{r_1}\cdots p_s^{r_s}p^r$ is the prime
factorization of $h$ and we have
$$\Phi_h(x)=\Phi_{p_1\cdots p_sp}(x^{(p_1^{(r_1-1)}\cdots p_s^{(r_s-1)}p^{(r-1)})}).$$
Put for short $y=x^{(p_1^{(r_1-1)}\cdots p_s^{(r_s-1)})}$, and
rewrite
$$\Phi_h(x)=\Phi_{p_1\cdots p_sp}(y^{p^{(r-1)}})=\frac{\Phi_{p_1\cdots p_s}(y^{p^r})}
{\Phi_{p_1\cdots p_s}(y^{p^{(r-1)}})}.$$ Since we are evaluating
on a field of characteristic $p$, we have now
$$\Phi_h(x)=\frac{\Phi_{p_1\cdots p_s}(y)^{p^r}}
{\Phi_{p_1\cdots p_s}(y)^{p^{(r-1)}}}=\Phi_{p_1\cdots
p_s}(y)^{p^{(r-1)}(p-1)}.$$ But we have, as well
$$\Phi_k(x)=\Phi_{p_1\cdots p_s}(x^{(p_1^{(r_1-1)}\cdots p_s^{(r_s-1)})})=\Phi_{p_1\cdots p_s}(y),$$
so we are done.

\end{document}